\documentstyle[11pt,amscd,amsfonts,leqno]{amsart}

\newtheorem{Theorem}{Theorem}[section]
\newtheorem{Proposition}[Theorem]{Proposition}

\newtheorem{Corollary}[Theorem]{Corollary}
\newtheorem{Definition}[Theorem]{Definition}
\newtheorem{Example}[Theorem]{Example}
\newtheorem{Remark}[Theorem]{Remark}

\begin{document}

\title{Quasiplurisubharmonic Green functions}
\author{Dan Coman and Vincent Guedj}
\date{}


\thanks{D. Coman was supported by the NSF grant DMS 0500563. Support by the Institut Mittag-Leffler (Djursholm, Sweden) is gratefully acknowledged.}

\address{D. Coman: dcoman@@syr.edu, Department of Mathematics, Syracuse
University, Syracuse, NY 13244-1150, USA}
\address{V. Guedj: guedj@@cmi.univ-mrs.fr, Universit\'e Aix-Marseille 1, LATP, 13453 Marseille Cedex 13, FRANCE}

\pagestyle{myheadings} \markboth{Dan Coman and Vincent Guedj}{Quasiplurisubharmonic Green functions}

\begin{abstract}
\noindent Given a compact K\"ahler manifold $X$, a quasiplurisubharmonic function is called a Green function with pole at $p\in X$ if its Monge-Amp\`ere measure is supported at $p$. We study in this paper the
existence and properties of such functions, in connection to their singularity at $p$. A full characterization is obtained in concrete cases, such as (multi)projective spaces.
\end{abstract}

\maketitle

\section*{Introduction}
Let $X$ be a compact K\"ahler manifold of complex dimension $n$. We pursue the study started in \cite{Yau}, \cite{Ko98}, \cite{Ko05}, \cite{GZ2}, \cite{EGZ}, \cite{BGZ} of the range of the complex Monge-Amp\`ere operator. 
Given a K\"ahler class $\alpha\in H^{1,1}(X,{\Bbb R})$ and a positive Radon measure $\mu$, the problem is to solve the equation $T^n=\mu$, where $T$ is a positive closed (1,1)-current in $\alpha$. When $\mu$ does not charge pluripolar sets, a complete answer was given in \cite{GZ2}. The main purpose of this article is to start and study the case when $\mu$ charges pluripolar sets by looking at measures $\mu$ which are sums of Dirac masses. The equation now reads \begin{equation}\label{e:MAeq}
T^n=\sum_{j=1}^k c_j \delta_{p_j}.
\end{equation}
We  seek solution(s) $T \in\alpha$ whose potentials are locally bounded away from the poles $p_1,\dots,p_k$. An obvious necessary condition in order to solve (\ref{e:MAeq}) is that the volume of $\alpha$,
$$V_\alpha:={\rm Vol}\,(\alpha)=\alpha^n,$$ 
is equal to the total mass of $\mu$, $\mu(X)=\sum c_j={\rm Vol}\,(\alpha)$.

\medskip

\par Fix $\theta$ a K\"ahler form representing $\alpha$ and let $PSH(X,\theta)$ denote the set of $\theta$-plurisubharmonic ($\theta$-psh) functions: these are functions $\varphi \in L^1(X,{\Bbb R})$ which are upper semicontinuous and such that $T=\theta+dd^c\varphi$ is a positive current. Here $d=\partial+\overline{\partial}$ and $d^c=\frac{1}{2\pi i} ( \partial-\overline{\partial})$. Solving (\ref{e:MAeq}) is therefore equivalent to finding a { \em``quasiplurisubharmonic Green function"}:

\vskip3mm

\noindent {\bf Definition.} {\em A function $\varphi \in PSH(X,\theta)$ is called a $\theta$-psh Green function with (isolated) poles at $p_1,\dots,p_k \in X$ if it is locally bounded in $X \setminus\{p_1,\ldots,p_k\}$ and  
$$(\theta+dd^c\varphi)^n =V_\alpha\sum_{j=1}^k m_j \delta_{p_j},\;\text{where}\;m_j>0,\;\sum_{j=1}^k m_j=1.$$}

\par In \cite{CGZ}, the domain $DMA(X,\theta)$ of the Monge-Amp\`ere operator was defined as the largest set of $\theta$-psh functions on which the operator is continuous along decreasing sequences of bounded 
$\theta$-psh functions. Hence one can consider a more general notion of $\theta$-psh Green function, by only requiring in the above definition that $\varphi\in DMA(X,\theta)$, instead of $\varphi$ being locally bounded away from the poles. We will not pursue this here. 

\medskip

\par Similar objects were considered by several authors in a local context  (\cite{Lem83}, \cite{Kli85}, \cite{D87}, \cite{Lel89}, \cite{CeP}, \cite{CoMZ}, \cite{CN}), and have found important applications (see e.g. \cite{BP98}, \cite{He99}, \cite{DiH00}). In our global context their existence depends on the geometry of $X$ and on the local positivity properties of $\alpha$ at the poles.

\medskip

\par We therefore study in {\em section \ref{S:locpos}} several indicators of the local positivity properties of 
$\alpha$, following Demailly \cite{D90B}. Recall that the Lelong number $\nu(\varphi,x)$ of a $\theta$-psh function $\varphi$ at $x$ is the largest constant $\nu$ for which $\varphi(p)\leq\nu\log dist(p,x)+O(1)$ holds for $p$ near $x$. If $\varphi(p)=\nu\log\,dist(p,x)+O(1)$ for $p$ near $x$ and $\nu>0$, we say that $\varphi$ has an {\em isotropic pole} at $x$ with Lelong number $\nu$.

\par We let $\nu(\alpha,x)$ (resp. $\varepsilon(\alpha,x)$) denote the maximal (resp. maximal isotropic) logarithmic singularity that a positive closed current $T \in \alpha$ can have at the point $x$. The indicator 
$\varepsilon(\alpha,x)$, introduced by Demailly \cite{D90B}, is called the Seshadri constant of $\alpha$ at 
$x$ and was intensively studied in algebraic geometry. We note in {\em section \ref{S:locpos}} that for all 
$x \in X$,
$$\nu(\alpha,x) \geq{\rm Vol}\,(\alpha)^{1/n} \geq \varepsilon(\alpha,x).$$

\par Thus a necessary condition for the existence of a $\alpha$-Green function with one  isotropic pole at $x$ is that ${\rm Vol}\,(\alpha)^{1/n} =\varepsilon(\alpha,x)$. This is far from being true in general: we observe for instance in Proposition \ref{P:LelP1} that this is never the case when $X$ is a multiprojective space. Even if this condition is satisfied, it is not clear whether it is sufficient, nor is it clear that the supremum in the definition of 
$\varepsilon$ is attained. We observe in section \ref{SSS:DPgen} that the following properties are equivalent:
\begin{itemize}
\item existence of a Green function with 9 isotropic poles in general position in ${\Bbb P}^2$;
\item existence of a Green function with one isotropic pole in generic position on a degree 1 Del Pezzo surface;
\item existence of a positive metric with bounded potentials for $c_1(Y)$, where $Y \rightarrow {\Bbb P}^2$ denotes the blow up of ${\Bbb P}^2$ at 9 points in general position,
\end{itemize}
the last one being a famous open problem \cite{DPS}. We therefore introduce in {\em section \ref{S:locpos}} weaker notions of Green functions. We show in Theorems \ref{T:Kahcon}, \ref{T:Kahcons} and Proposition \ref{P:pLelX} how to construct these by a balayage procedure. It is a delicate and interesting problem to determine whether $\theta$-psh Green functions always exist. As already observed, we have to consider arbitrary singularities. The balayage procedure depends on the choice of local data $(u_1,\dots,u_k)$ encoding the singularities at the poles $(p_1,\dots,p_k)$. In particular, the problem of constructing $\theta$-psh Green functions is reduced to finding local data for which the functions $g$ constructed in Theorems \ref{T:Kahcon} and \ref{T:Kahcons} have isolated singularities at $p_j$. 

\medskip

\par In {\em section \ref{S:Pn}} we give a complete description of all these notions on the complex projective space ${\Bbb P}^n$. In particular, we characterize in Theorem \ref{T:rigid} Green functions arising naturally from rational maps $f:{\Bbb P}^n\dashrightarrow{\Bbb P}^{n-1}$ with finite indeterminacy set. We end section \ref{S:Pn} by constructing interesting dynamical Green functions.

\medskip

\par In {\em section \ref{S:P1P1}} we compute similar quantities for multiprojective spaces, focusing on 
${\Bbb P}^1 \times{\Bbb P}^1$. We show in Proposition \ref{P:GreenP1} that Green functions with one pole correspond to a certain class of Green functions with three poles on ${\Bbb P}^2$. A large class of examples of these can be constructed using Theorem \ref{T:rigid} (see Example \ref{E:GreenP1}). However, there is no Green function with one isotropic pole on ${\Bbb P}^1 \times{\Bbb P}^1$ (Corollary \ref{C:LelP1}). 

\medskip
 
\par In {\em section \ref{S:DelPezzo}} we turn our attention to the case of smooth Del Pezzo surfaces, focusing on those of degree 1, i.e. blow ups $X$ of ${\Bbb P}^2$ at 8 points in general position. Let $\alpha$ be  the first Chern class of $X$. We prove in Proposition \ref{P:LX8} that $\nu(\alpha,x)=1$ if $x\in X\setminus S$, and $\nu(\alpha,x)=2$ if $x\in S$. Here $S$ is the set of singular points on the singular cubics passing through the 8 blown up points, and $1\leq|S|\leq12$. The results of Proposition \ref{P:LX8} allow us to compute, using currents, the exact value of Tian's ``$\alpha$-invariant", and to deduce that $X$ has a K\"ahler-Einstein metric (section \ref{SS:DPKE}). We conclude the paper with the discussion in section \ref{SS:DPGreen} of $\omega$-psh Green functions with one pole $x\in X$, where $\omega\in\alpha$ is a 
K\"ahler form. Such functions are easy to construct when $x\in S$. For generic points $x\not\in S$ the existence of Green functions with an isotropic pole at $x$ of maximal Lelong number 
$1=\varepsilon(\alpha,x)$ is equivalent to a famous open problem in algebraic geometry (see section {\ref{SSS:DPgen}).

\medskip

\noindent {\bf Acknowledgement.} We would like to thank the referee for his comments which helped improve the exposition of this paper.

\section{Local positivity of (1,1) classes and Green functions}\label{S:locpos}
Let ${\mathcal P}(X)$ be the set of all positive closed currents of bidegree (1,1) on $X$. For $\alpha\in H^{1,1}(X,{\Bbb R})$ we let $${\mathcal P}(\alpha)=\{T\in{\mathcal P}(X):\,T\in\alpha\}$$ be the set of positive closed currents whose cohomology class is $\alpha$. By definition, a class $\alpha$ is {\em pseudoeffective}  if ${\mathcal P}(\alpha)\neq\emptyset$. Let $H^{1,1}_{psef}(X,{\Bbb R})$ denote the closed convex cone of all pseudoeffective (1,1) classes. 

\smallskip

\par There are two other interesting cones in $H^{1,1}_{psef}(X,{\Bbb R})$ which correspond to stronger notions of positivity. We let $H^{1,1}_{Kaehler}(X,{\Bbb R})$ denote the cone of K\"ahler classes and  $H^{1,1}_{nef}(X,{\Bbb R})$ denote its closure. Then 
$H^{1,1}_{Kaehler}(X,{\Bbb R})$ is the interior of $H^{1,1}_{nef}(X,{\Bbb R})$.

\par Following Demailly \cite{D90B}, we would like to measure the local positivity of a class $\alpha$. There are two main indicators, in connection to the various types of positivity. In the sequel we denote by $\nu(T,x)$ the Lelong number of $T\in{\mathcal P}(X)$ at a point $x$.

\begin{Definition}\label{D:lpd} Let $\pi:\widetilde X\to X$ denote the blow up of $X$ at a point $x$, and let $E=\pi^{-1}(x)$ denote the exceptional divisor. 

\par 1) For $\alpha\in H^{1,1}_{psef}(X,{\Bbb R})$ we set 
$$\nu(\alpha,x):=\sup\{\nu\geq0:\,\pi^\star\alpha-\nu E\in H^{1,1}_{psef}(\widetilde X,{\Bbb R})\}.$$

\par 2) For $\alpha\in H^{1,1}_{nef}(X,{\Bbb R})$ we set 
$$\varepsilon(\alpha,x):=\sup\{\varepsilon\geq0:\,\pi^\star\alpha-\varepsilon E\in H^{1,1}_{nef}(\widetilde X,{\Bbb R})\}.$$
\end{Definition}

\par The indicator $\nu(\alpha,x)$ is the maximal Lelong number that a current $T\in{\mathcal P}(\alpha)$ can have at $x$. In this case the supremum is attained, because ${\mathcal P}(\alpha)$ is a compact set (in the weak topology of currents). 

\par The indicator $\varepsilon(\alpha,x)$ is called the Seshadri constant of $\alpha$ at $x$. It has been intensively studied since it was introduced by Demailly. We refer the reader to \cite[Chapter 5]{La} for a detailed account of this notion.

\smallskip

\par By definition we have $0\leq\varepsilon(\alpha,x)\leq\nu(\alpha,x)$. It follows from the characterization of the K\"ahler cone obtained in \cite{DP04} that if $\alpha\in H^{1,1}_{nef}(X,{\Bbb R})$ and $x\in X$ then 
$$\varepsilon(\alpha,x)=\min_V\left(\frac{(\alpha^{\dim V}\cdot V)}{{\rm mult}_x\,V}\right)^{\frac{1}{\dim V}},$$
where the minimum is taken over all irreducible subvarieties $V\subseteq X$ with $\dim V\geq1$ and $x\in V$ (see e.g. Proposition 5.1.9 and Remark 1.5.32 in \cite{La}). With $V=X$, this yields the estimate (recall that $V_\alpha={\rm Vol}\,(\alpha)$):
\begin{equation}\label{e:Seshvar}
\varepsilon(\alpha,x)\leq V_\alpha^{1/n},\;\forall\,x\in X.
\end{equation}
On the other hand, it follows easily from Theorem \ref{T:Kahcon} below that if  
$\alpha\in H^{1,1}_{Kaehler}(X,{\Bbb R})$  
$$ \nu(\alpha,x)\geq V_\alpha^{1/n},\;\forall\,x\in X.$$
Both bounds are sharp in the case of ${\Bbb P}^n$. 

\begin{Remark}\label{R:vamp} If $\alpha\in H^2(X,{\Bbb Z})$ is an integral class, then $\nu(\alpha,x)\geq V_\alpha^{1/n}\geq1$ for all $x\in X$. Note also that if $\alpha$ is very ample then $\varepsilon(\alpha,x)\geq1$. 
\end{Remark}

\par An alternate description of the Seshadri constant $\varepsilon(\alpha,x)$ can be given in terms the maximal Lelong number of currents in ${\mathcal P}(\alpha)$ whose potentials have an isolated singularity at 
$x$ \cite{D90B}. Let $\alpha\in H^{1,1}_{Kaehler}(X,{\Bbb R})$ and $\theta$ be a K\"ahler form representing $\alpha$. It follows as in \cite[Theorem 6.4]{D90B} that for every $x\in X$, 
\begin{eqnarray}\label{e:Seshadri}
\varepsilon(\alpha,x) & = & \sup\{\gamma:\,\exists\,\varphi\in PSH(X,\theta),\;\|\varphi-\gamma\log dist(\cdot,x)\|_{L^\infty(X)}<+\infty\}\\
& = & \sup\{\gamma:\,\exists\,\varphi\in PSH(X,\theta),\;\nu(\varphi,x)=\gamma,\;\varphi\in L^\infty_{loc}(U\setminus\{x\})\},\nonumber
\end{eqnarray}
where $U$ is a neighborhood of $x$ depending on $\varphi$. Recall that $PSH(X,\theta)$ is the set of $\theta$-psh functions. The set of normalized $\theta$-psh functions, for example by the condition $\max_X\varphi=0$, is isomorphic to ${\mathcal P}(\alpha)$ via $\varphi\to\theta+dd^c\varphi\in{\mathcal P}(\alpha)$. The fact that the two supremums are equal is straightforward. Moreover, in this case we have $\varepsilon(\alpha,x)>0$ for all $x\in X$.

\medskip

\par We now list  a few elementary properties of these numerical indicators.  

\begin{Proposition}\label{P:lpp} 1) The functions $\alpha\to\nu(\alpha,x),\varepsilon(\alpha,x)$ are homogeneous and superadditive (i.e. $\nu(\alpha+\beta,x)\geq\nu(\alpha,x)+\nu(\beta,x)$).

\par 2) The function $x\to\nu(\alpha,x)$ is upper semicontinuous. 
 
\par 3) If $\alpha$ is K\"ahler the function $x\to\varepsilon(\alpha,x)$ is lower semicontinuous.
\end{Proposition}

\begin{pf} The upper semicontinuity property of $x\to\nu(\alpha,x)$ follows since ${\mathcal P}(\alpha)$ is compact and from the well known fact that $\limsup\nu(T_j,x_j)\leq\nu(T,x)$ as positive closed (1,1)-currents $T_j\to T$ and $x_j\to x$.

\par To prove 3), let $\theta\in\alpha$ be a K\"ahler form, $x\in X$, $0<\epsilon<1$, and $0<\nu<\varepsilon(\alpha,x)$. We construct for all $y$ near $x$ a $\theta$-psh function
$\varphi_y$ with
$\varphi_y=(1-\epsilon)\nu\log dist(\cdot,y)+O(1)$.
Using (\ref{e:Seshadri}), this shows that 
$\liminf_{y\to x}\varepsilon(\alpha,y)\geq\varepsilon(\alpha,x)$.

\par By (\ref{e:Seshadri}) there exists $\varphi\in PSH(X,\theta)$ such that
$\varphi=\nu\log dist(\cdot,x)+O(1)$. Let $B_2\subset{\Bbb C}^n$ be the ball of radius 2 centered at 0. We can find a coordinate chart $f:B_2\longrightarrow U\subset X$, $f(0)=x$, and a function $\rho\in C^\infty(U)$ so that $dd^c\rho=\theta$ and 
$$\nu\log\|z\|-C\leq v(z):=(\rho+\varphi)\circ f(z)
\leq\nu\log\|z\|+C,\;z\in B_2,$$ for some constant $C>0$. Fix $r>0$
small enough so that
$$(1-\epsilon)\left(\nu\log\frac{r}{2}-2C\right)\geq
\nu\log r+2C.$$ Next, let $T_w$ be an automorphism of the unit
ball $B_1\subset{\Bbb C}^n$ with $T_w(w)=0$. There exists $\delta(r)<r$ such
that $\|T_w(z)\|\geq r/2$, if $\|z\|=r$ and $\|w\|<\delta(r)$. For
such $w$ we define the function $v_w$ on $B_2$ by
$$v_w(z)=\left\{ \begin{array}{ll}v(z)+C,\;1\leq\|z\|<2,\\
\max\{v(z)+C,(1-\epsilon)(v\circ T_w(z)-C)\},\;r<\|z\|<1,\\
(1-\epsilon)(v\circ T_w(z)-C),\;\|z\|\leq r.
\end{array}\right.$$ Note that if $\|z\|=1$ then
$v(z)+C\geq0\geq(1-\epsilon)(v\circ T_w(z)-C)$, while if
$\|z\|=r$, $$(1-\epsilon)(v\circ T_w(z)-C)\geq(1-\epsilon)
\left(\nu\log\frac{r}{2}-2C\right)\geq\nu\log r+2C\geq v(z)+C.$$
Hence $v_w$ is psh on $B_2$ and
$v(z)=(1-\epsilon)\nu\log\|z-w\|+O(1)$ for $z$ near $w$.

\par For $y=f(w)$, where $\|w\|<\delta(r)$, we finally let
$$\varphi_y=\left\{\begin{array}{ll}\varphi+C,
\;{\rm on}\;X\setminus f(B_1),\\v_w\circ f^{-1}-\rho,\;{\rm
on}\;f(B_1).\end{array}\right.$$ Then $\varphi_y$ is
$\theta$-psh and
$\varphi_y=(1-\epsilon)\nu\log dist(\cdot,y)+O(1)$ near $y$.
\end{pf}

\par In general, the functions $\nu(\alpha,\cdot),\varepsilon(\alpha,\cdot)$ are not continuous (see e.g. Proposition \ref{P:LX8} and section \ref{SS:DPGreen}). Note that in the special case when $X$ is projective and $\alpha$ is an integral class, it follows from \cite[Example 5.1.11]{La} that $\varepsilon(\alpha,\cdot)$ is constant outside a countable union of proper subvarieties of $X$. 

\smallskip

\par If $\theta\in\alpha$ is a K\"ahler form, we have by (\ref{e:Seshvar}) and (\ref{e:Seshadri}) that a necessary condition for the existence of a $\theta$-psh Green function with an {\em isotropic} pole at $p$ is 
$$\varepsilon(\alpha,p)=V_\alpha^{1/n}.$$ 
Since this fails to hold in general (see Proposition \ref{P:LelP1}), one has to consider other singularities. Following ideas of Demailly \cite{D93b}, we will show that local fundamental solutions of the Monge-Amp\`ere operator have $\theta$-psh subextensions to $X$. 

\smallskip

\par We will consider the slightly more general situation when the class $\alpha$ is represented by a smooth closed (1,1) form $\theta\geq0$ and $V_\alpha>0$. Recall that the {\em unbounded locus} $M(\varphi)$ of $\varphi\in PSH(X,\theta)$ is defined as the set of all points
$p\in X$ such that $\varphi$ is unbounded in every neighborhood of
$p$. We denote by $PSH^-(X,\theta)$ the set of $\theta$-psh functions $\varphi\leq0$ on $X$. For $p\in X$, let 
${\mathcal G}_p(V_\alpha)$ be the set of germs of functions $u$ at $p$ with the
following properties: there exists an open set $U\subset X$
containing $p$ such that $u$ is psh on $U$ and locally bounded  on
$U\setminus\{p\}$, $u(p)=-\infty$, and
$(dd^cu)^n=V_\alpha\delta_p$ as measures on $U$. 

\begin{Theorem}\label{T:Kahcon} Let $p\in X$ and $u\in{\mathcal
G}_p(V_\alpha)$. There exists a unique function $g=g_{u,p}\in
PSH^-(X,\theta)$ such that

\par (i) $g\leq u+C$ holds near $p$, for some constant C.

\par (ii) If $\varphi\in PSH^-(X,\theta)$ and $\liminf_{q\to
p}\varphi(q)/u(q)\geq 1$ then $\varphi\leq g$ on $X$.
\\ In addition, $g$ has the following properties:

\par (a) $(\theta+dd^cg)^n=0$ on the open set
$X\setminus(M(g)\cup\{g=0\})$.

\par (b) If $p$ is an isolated point of $M(g)$ then $M(g)=\{p\}$
and $g$ is a $\theta$-psh Green function on $X$ with pole at $p$.

\par (c) The open set $D_{u,p}=\{g<0\}$ is connected.
\end{Theorem}

\par It should be noted that the existence of a global $\theta$-psh function
$\varphi$ subextending $u$ (i.e. such that $\varphi \leq u$ near $p$) is a nontrivial 
matter. We use Yau's solution in the spirit of \cite{D93b}, \cite{DP04}.
Producing the ``best subextension'' $g$ proceeds using a classical balayage procedure 
(see \cite{Ra06} for recent similar local extremal problems).

\begin{pf} The uniqueness of a function with properties $(i),\,(ii)$
is clear. Fix $U\subset X$ an open coordinate ball around $p$, so
that $u$ is psh on $U$, locally bounded on $U\setminus\{p\}$ and
$(dd^cu)^n=V_\alpha\delta_p$ as measures on $U$. We divide the
proof in three steps.

\par{\bf Step 1.} Using a mass concentration technique of Demailly \cite{D93b}, we construct a function
$\varphi\in PSH(X,\theta)$ so that $\varphi\leq u$ near $p$. Let $\omega_0$ be a K\"ahler form on $X$.

\par Let $W\subset\subset W'\subset\subset U$ be open and connected, with
$p\in W$, and let $\chi$ be a smooth function on $X$ with compact
support in $W'$, such that $0\leq\chi\leq1$ and $\chi=1$ on $W$.
We may assume that $u\geq0$ on $\partial W$. Let $\rho,\,\rho_0$ be negative smooth functions on 
$W'$ with $dd^c\rho=\theta$, $dd^c\rho_0=\omega_0$.

\par Let $u_j\searrow u$ be a sequence of smooth psh functions 
on $W'$ and let $\omega_j=\theta+j^{-1}\omega_0$. We define measures
$$\mu_j=C_j\chi\,(dd^c\,u_j)^n,$$ where the constants $C_j>0$ are 
chosen so that $\mu_j(X)=\int_X\omega_j^n$. Note that $\mu_j$ has
support in $W'$, and $(dd^c\,u_j)^n\to V_\alpha\delta_p$ in the
weak sense of measures on $W'$. Hence
$$\lim_{j\rightarrow\infty}\int\chi\,(dd^c\,u_j)^n=V_\alpha\chi(p)
=V_\alpha,\;{\rm so}\;\lim_{j\rightarrow\infty}C_j=1.$$

\par Yau's theorem (see \cite{Yau}, also \cite{Ko98}) implies that
there exist continuous functions $\varphi_j\in PSH(X,\omega_j)$
such that
$$(\omega_j+dd^c\varphi_j)^n=\mu_j,\;\max_X\,\varphi_j=0.$$ 
By \cite[Proposition 1.7]{GZ1} we may assume after passing to a subsequence that $\{\varphi_j\}$ converges in $L^1(X)$ to a function 
$\varphi\in PSH(X,\theta)$. Moreover, by \cite[Theorem 4.1.8]{Ho}
we have $\varphi=(\limsup_{j\rightarrow\infty}\varphi_j)^\star$ on $X$.

\par Choose a sequence $a_j\geq1$ so that
$a_j^nC_j>1$ and $a_j\to1$. We have
$$a_j(\varphi_j+\rho+j^{-1}\rho_0)\leq0
\leq u_j\;{\rm on}\;\partial W.$$ On the other hand
$$a_j^n(dd^c(\varphi_j+\rho+j^{-1}\rho_0))^n=
a_j^nC_j\chi\,(dd^c\,u_j)^n\geq(dd^c\,u_j)^n$$ holds on $W$, as
$\chi=1$ on $W$. The minimum principle of Bedford and Taylor
\cite[Theorem A]{BT76} implies that
$a_j(\varphi_j+\rho+j^{-1}\rho_0)\leq u_j$ on $W$. Letting
$j\rightarrow\infty$ we obtain that $\varphi+\rho\leq u$ holds on
$W$.  This concludes Step 1.

\par {\bf Step 2.} We construct the function $g$ using an
upper envelope method. Consider the family
$${\mathcal F}=\left\{\varphi\in PSH^-(X,\theta):\,\liminf_{q\to
p}\frac{\varphi(q)}{u(q)}\geq 1\right\}.$$
In the terminology of Rashkovskii, this is the family of negative $\theta$-psh functions whose relative type with respect to $u$ is at least 1 (see \cite{Ra06}).

By Step 1, ${\mathcal F}\neq\emptyset$. If $g=\sup\{\varphi:\,\varphi\in{\mathcal F}\}$, then the upper
semicontinuous regularization $g^\star\in PSH^-(X,\theta)$. We will show that $g^\star\leq u+C$ holds near $p$ for some constant $C$. This implies that $g=g^\star\in{\mathcal F}$, so $g$ verifies properties $(i),\,(ii)$.

\par We can find $M>0$ such that the connected component $D$ of
$\{u<-M\}$ which contains $p$ is relatively compact in $U$. Let
$\rho<0$ be a smooth function on $U$ so that $dd^c\rho=\theta$. Fix $\varphi\in{\mathcal F}$. There exists a sequence of
relatively compact domains $D_j\subset D$, $j>0$, with the
following properties: $$D_{j+1}\subset
D_j,\;\bigcap_{j>0}D_j=\{p\},\; \varphi(q)\leq
(1-j^{-1})u(q)\;\text{for }q\in\overline D_j.$$ We have
$\rho+\varphi\leq0\leq (1-j^{-1})(u+M)$ on $\partial D$, and
clearly $\rho+\varphi\leq(1-j^{-1})(u+M)$ on $\partial D_j$. Since
the psh function $u$ is maximal on $U\setminus\{p\}$, it follows
that the last inequality holds on $D\setminus D_j$. As
$j\to\infty$ we see that $\rho+\varphi\leq u+M$ on $D$. Since
$\varphi\in{\mathcal F}$ was arbitrary, this implies that
$g^\star\leq u+C$ on $D$, where $C=M-\min_D\rho$.

\par {\bf Step 3.} We prove the remaining properties of $g$.

\par $(a)$ Note that $M(g)$ is closed and since $g\leq0$ is upper semicontinuous the set
$\{g=0\}$ is closed. Let $q\in X\setminus(M(g)\cup\{g=0\})$ and
let $\rho$ be a smooth function in a neighborhood
of $q$ such that $dd^c\rho=\theta$ and $\rho(q)=0$. We can find $\varepsilon>0$ and a
small neighborhood $G$ of $q$ such that $G\subset
X\setminus(M(g)\cup\{g=0\})$ and $g<-\varepsilon$,
$|\rho|<\varepsilon/2$ on $G$. Let $W$ be a relatively compact
open subset of $G$ and $v$ be psh on $W$ so that
$v^\star\leq\rho+g$ on $\partial W$. The function
$$\varphi=g\;\;\text{on}\;X\setminus W,\;\;\;
\varphi=\max\{\rho+g,v\}-\rho\;\;\text{on}\;W,$$ is $\theta$-psh
and $\varphi\leq0$ on $X$. Since $\varphi=g$ in a neighborhood of
$p$, we conclude that $\varphi\in{\mathcal F}$, hence $v\leq\rho+g$ on $W$. This shows that the psh function
$\rho+g$ is maximal on $G$. By \cite{BT82}, $(\theta+dd^cg)^n=0$ in
$G$, and hence on $X\setminus(M(g)\cup\{g=0\})$.

\par $(b)$ If $p\in M(g)$ is isolated, there exists
a closed ball $K$ centered at $p$ so that $K\cap M(g)=\{p\}$.
Hence $g$ is bounded below on $\partial K$. It follows that if $C>0$ is large enough the
function $\varphi$ defined by $\varphi=g$ on $K$,
$\varphi=\max\{g,-C\}$ on $X\setminus K$, is $\theta$-psh and
$\varphi\in{\mathcal F}$. Thus $\varphi\leq g$, so $M(g)=\{p\}$.
By $(i)$ and \cite{D93},
$(\theta+dd^cg)^n(\{p\})\geq(dd^cu)^n(\{p\})=V_\alpha$. Mass considerations imply that $g$ is a
$\theta$-psh Green function.

\par $(c)$ Suppose that there exists a connected component $W$ of
$D_{u,p}$ not containing $p$. The function $\varphi$ defined by
$\varphi=g$ on $X\setminus W$ and $\varphi=0$ on $W$, verifies
$\varphi\in{\mathcal F}$, so $\varphi\leq g$. This contradicts our
assumption that $g<0$ on $W$, so $D_{u,p}$ is
connected.
\end{pf}

\par The following theorem produces Green functions with several poles. Its proof is a straightforward adaptation of the proof of Theorem \ref{T:Kahcon}.

\begin{Theorem}\label{T:Kahcons} For $1\leq j\leq k$, let $p_j\in X$, $u_j\in{\mathcal G}_{p_j}(V_\alpha)$, and $m_j>0$ with $\sum_{j=1}^km_j=1$. There exists a unique function $g\in PSH^-(X,\theta)$ such that

\par (i) $g\leq m_j^{1/n}u_j+C$ holds near each $p_j$, for some constant C.

\par (ii) If $\varphi\in PSH^-(X,\theta)$ and for each $j$, $\liminf_{q\to
p_j}\varphi(q)/u_j(q)\geq m_j^{1/n}$, then $\varphi\leq g$ on $X$.

\par Moreover, we have $(\theta+dd^cg)^n=0$ on $X\setminus(M(g)\cup\{g=0\})$.
If all $p_j$ are isolated points of $M(g)$ then $g$ is a $\theta$-psh Green function with poles at 
$p_1,\dots,p_k$.
\end{Theorem}

\medskip

\par It is an intricate problem to decide whether there always exist  local models $u$ at $p\in X$ such that $g_{u,p}$ is a Green function. As an alternate approach, we introduce a partial Green function associated to an isotropic singularity.  

\begin{Proposition}\label{P:pLelX} Let  $\theta\in\alpha$ be a K\"ahler 
form, let $p\in X$ and $0<\gamma<\varepsilon(\alpha,p)$. There exists a unique function
$\psi_{\gamma,p}\in PSH^-(X,\theta)$ so that
$\nu(\psi_{\gamma,p},p)=\gamma$ and with the property that if
$\varphi\in PSH^-(X,\theta)$ and $\nu(\varphi,p)\geq\gamma$ then
$\varphi\leq\psi_{\gamma,p}$. Moreover,  
$$\|\psi_{\gamma,p}-\gamma\log dist(\cdot,p)\|_{L^\infty(X)}<+\infty,\;
(\theta+dd^c\psi_{\gamma,p})^n=\gamma^n\delta_p+\mu_{\gamma,p},$$ 
where $\mu_{\gamma,p}$ is a positive measure supported on the compact 
$\{\psi_{\gamma,p}=0\}$.
\end{Proposition}

\begin{pf} The uniqueness of $\psi_{\gamma,p}$ is clear. Let us fix a
biholomorphic map $f:B\to U$
from the unit ball $B\subset{\Bbb C}^n$ onto a neighborhood $U$ of $p$, with $f(0)=p$. Let $\rho<0$ 
be a smooth function on $U$ with $dd^c\rho=\theta$.

\par By (\ref{e:Seshadri}) there exists $\psi\in PSH^-(X,\theta)$ so that $\psi=\gamma\log dist(\cdot,p)+O(1)$. Let 
$$\psi_{\gamma,p}(q)=\sup\{\varphi(q):\,\varphi\in PSH^-(X,\theta),\;\nu(\varphi,p)\geq\gamma\}.$$ 
For such $\varphi$, we have $(\rho+\varphi)(f(z))\leq\gamma\log\|z\|$ on $B$. This implies 
$\psi_{\gamma,p}^\star\in PSH^-(X,\theta)$ and
$\nu(\psi_{\gamma,p}^\star,p)\geq\gamma$. Thus
$\psi_{\gamma,p}=\psi_{\gamma,p}^\star$. Since $\psi\leq\psi_{\gamma,p}$, it follows that 
$\nu(\psi_{\gamma,p},p)=\gamma$ and the function 
$\psi_{\gamma,p}-\gamma\log dist(\cdot,p)$ is bounded on $X$.

\par Arguing as in the proof of Theorem \ref{T:Kahcon} $(a)$ we show that 
$(\theta+dd^c\psi_{\gamma,p})^n=0$ in $\{\psi_{\gamma,p}<0\}\setminus\{p\}$. By \cite{D93}, 
$(\theta+dd^c\psi_{\gamma,p})^n(\{p\})=\gamma^n$, and the proof is complete.
\end{pf}

\par We refer to \cite{Ra06} for similar extremal problems on domains in ${\Bbb C}^n$. In the following sections, we are going to compute the functions $\nu,\,\varepsilon$ and $g_{u,p},\,\psi_{\nu,p}$ in a number of interesting cases.

\section{Green functions on ${\Bbb P}^n$}\label{S:Pn}
\par Let $[z_0:\ldots:z_n]$ be homogeneous coordinates on
${\Bbb P}^n$ and $\pi_n:{\Bbb
C}^{n+1}\setminus\{0\}\rightarrow{\Bbb P}^n$ be the standard
projection. Let $\alpha_n=\{\omega_n\}$, where $\omega_n$ is the Fubini-Study
form, so $\pi_n^\star\omega_n=dd^c\log\|z\|$ and ${\rm Vol}\,(\alpha_n)=1$. 

\subsection{Maximal Lelong number}\label{SS:PnMax}

\begin{Proposition}\label{P:mLPn} We have $\nu(\alpha_n,x)=\varepsilon(\alpha_n,x)=1$ for all $x\in{\Bbb P}^n$. If $T\in{\mathcal P}(\alpha_n)$ and $\nu(T,x)=1$ then $T=\wp_x^\star S$, where $\wp_x:{\Bbb P}^n\dashrightarrow{\Bbb P}^{n-1}$ is the projection with center $x$ onto a hyperplane ${\Bbb P}^{n-1}\not\ni x$ and $S\in{\mathcal P}(\alpha_{n-1})$. Moreover, the following are equivalent:

\par (i) the potentials of $T$ have isotropic pole at $x$ with Lelong number 1.

\par (ii) $T$ has locally bounded potentials on ${\Bbb P}^n\setminus\{x\}$.

\par (iii) $S$ has bounded potentials. 
\end{Proposition}

\begin{pf} Let $\pi:X\to{\Bbb P}^n$ denote the blow up of ${\Bbb P}^n$ at $x$, and let $E$ be the exceptional divisor. The map $\Phi=\wp_x\circ\pi:X\to{\Bbb P}^{n-1}$ is a holomorphic fibration, whose fibers are the projective lines through $x$. Moreover, $\pi^\star\alpha_n-E=\Phi^\star\alpha_{n-1}$. 

\par If $\nu(T,x)=1$ then $\widetilde T=\pi^\star T-[E]$ is a positive closed (1,1)-current on $X$ in the cohomology class $\Phi^\star\alpha_{n-1}$. It follows that $\widetilde T=\Phi^\star S$ for some $S\in{\mathcal P}(\alpha_{n-1})$, hence $T=\wp_x^\star S$. The potentials of $T$ have isotropic pole at $x$ with Lelong number 1 if and only if $\widetilde T$ has bounded potentials, hence  if and only if $S$ has bounded potentials. 

\par It is well known that currents in ${\mathcal P}(\alpha_n)$ have Lelong number at most 1 at each point $x$. The above construction shows that $\nu(\alpha_n,x)=\varepsilon(\alpha_n,x)=1$.
\end{pf}

\par We now explore further the geometry of sublevel sets of high Lelong numbers, in the spirit of \cite{CoPAMS}. For $c>0$ and $T\in{\mathcal P}(\alpha_n)$ a theorem of Siu \cite{Siu74} states that 
$$E_c(T):=\{x\in{\Bbb P}^n:\,\nu(T,x)\geq c\}$$ 
is an algebraic subset of dimension at most $n-1$. We also consider the set
$$E^+_c(T):=\{x\in{\Bbb P}^n:\,\nu(T,x)>c\}.$$

\begin{Proposition}\label{P:E+c} The set $E^+_{n/(n+1)}(T)$ is contained in a 
hyperplane of ${\Bbb P}^n$.\end{Proposition}

\begin{pf} Let $T=\omega_n+dd^c\varphi$ and set $E_c(\varphi)=E_c(T)$ and $E^+_c(\varphi)=E^+_c(T)$. The proof is by induction on $n$. If $n=1$, $T$ is a probability measure, $\nu(T,p)=T(\{p\})$, so $E^+_{1/2}(T)$ contains at most one point.

\par Let $c_n=n/(n+1)$. If $n\geq2$ we assume for a contradiction
that $E^+_{c_n}(\varphi)$ contains the points $q,p_1,\dots,p_n$
in general position. Let $H$ be the hyperplane determined by
$p_1,\dots,p_n$, so $q\not\in H$. By a theorem of Siu
\cite{Siu74}, $T=c[H]+R$, where $0\leq c\leq1$ and $R\in{\mathcal P}((1-c)\alpha_n)$ has generic
Lelong number 0 along $H$. Thus
$$c_n<\nu(\varphi,q)=\nu(R,q)\leq1-c,\;\nu(R,p_j)=\nu(\varphi,p_j)-c>c_n-c,\;1\leq j\leq n.$$
Consider the current $S=R/(1-c)=\omega_n+dd^c\psi\in{\mathcal P}(\alpha_n)$. Since 
$c<1-c_n$,  
$$\nu(\psi,p_j)>\frac{c_n-c}{1-c}>\frac{2c_n-1}{c_n}=c_{n-1},\;1\leq j\leq n.$$ 
By \cite[Proposition 3.7]{D92}, there exist $\epsilon_k\searrow0$ and currents
$S_k=(1+\epsilon_k)\omega_n+dd^c\psi_k\geq0$, where $\psi_k$ have analytic singularities, such that $S_k\rightarrow S$ and 
$0\leq\nu(\psi,p)-\nu(\psi_k,p)\leq\epsilon_k$ for all $p\in{\Bbb P}^n$. Since $S$ does not charge $H$, it follows that $\psi_k\not\equiv-\infty$ on $H\equiv{\Bbb P}^{n-1}$. Hence 
$\psi_k\mid_{_H}\in PSH({\Bbb P}^{n-1},\omega_{n-1})$ and
$$\nu(\psi_k\mid_{_H},p_j)\geq\nu(\psi_k,p_j)>c_{n-1},\;1\leq j\leq n,$$ for $k$ sufficiently
large. This yields a contradiction, since by our induction hypothesis the set $E^+_{(n-1)/n}(\psi_k\mid_{_H})$ is contained in a hyperplane of ${\Bbb P}^{n-1}$.\end{pf}

\par The value $n/(n+1)$ in the previous theorem is sharp. Indeed,
let $S$ be a set of $n+1$ points $p_j\in{\Bbb P}^n$ in general
position, and let $[H_j]$ be the current of integration along the
hyperplane $H_j$ determined by $S\setminus\{p_j\}$. If
$T=([H_1]+\ldots+[H_{n+1}])/(n+1)$ then the set $E_{n/(n+1)}(T)=S$
is not contained in any hyperplane.

\medskip

\par We are now in position to make the result of Proposition \ref{P:mLPn} more precise, by giving a characterization of the currents $T$ for which $E_1(T)\neq\emptyset$.

\begin{Proposition}\label{P:E1} If $T\in{\mathcal P}(\alpha_n)$ and  
$E_1(T)\neq\emptyset$ then $E_1(T)$ is a $k$-dimensional linear subspace of ${\Bbb
P}^n$ for some integer $0\leq k\leq n-1$. Let $\wp$ denote the projection with center $E_1(T)$ 
onto a linear subspace $L\equiv{\Bbb P}^{n-k-1}$ such that $L\cap E_1(T)=\emptyset$. Then $T=\wp^\star S$ 
for a unique current $S\in{\mathcal P}(\alpha_{n-k-1})$, and $E_1(S)=\emptyset$.
\end{Proposition}

\begin{pf} Let $T=\omega_n+dd^c\varphi$ and $k\geq0$ be the largest integer for which there exist $k+1$ points $p_0,\dots,p_k\in E_1(T)$ in general position (i.e. not contained in a $(k-1)$-dimensional subspace).  Proposition \ref{P:E+c} implies $k\leq n-1$. Using an automorphism of ${\Bbb P}^n$, we may assume 
$p_0=[1:0:\ldots:0]$, $p_1=[0:1:\ldots:0]$, and so on. Consider
the projection $f_0$ of ${\Bbb P}^n$ with center $p_0$ onto the
hyperplane ${\Bbb P}^{n-1}\equiv\{z_0=0\}$. Proposition 
\ref{P:mLPn} shows that $\varphi=u+h_0\circ f_0$, where $h_0\in PSH({\Bbb
P}^{n-1},\omega_{n-1})$ and 
$$u([z_0:\ldots:z_n])=\frac{1}{2}\,
\log\frac{|z_1|^2+\ldots+|z_n|^2}{|z_0|^2+\ldots+|z_n|^2}\;.$$
It follows that $f_0(p_j)\in E_1(h_0)$,
$j=1,\dots,k$, and Proposition \ref{P:mLPn} can be applied to
$h_0$ and the point $f_0(p_1)$. Continuing like this we get
$$\varphi([z_0:\ldots:z_n])=\frac{1}{2}\,
\log\frac{|z_{k+1}|^2+\ldots+|z_n|^2}{|z_{0}|^2+\ldots+|z_n|^2}+
h([z_{k+1}:\ldots:z_n]),$$ with $h\in PSH({\Bbb
P}^{n-k-1},\omega_{n-k-1})$. The definition of $k$ implies $E_1(h)=\emptyset$, so  
$E_1(\varphi)=\{z_{k+1}=\ldots=z_n=0\}$.\end{pf}

\subsection{Green functions}\label{SS:GPn}
\subsubsection{Green functions with one pole}
It follows from Proposition \ref{P:mLPn} that if $T=\wp_x^\star S$, where $S\in{\mathcal P}(\alpha_{n-1})$ has bounded potentials and $\wp_x:{\Bbb P}^n\dashrightarrow{\Bbb P}^{n-1}$ is the projection from $x$, then $T=\omega_n+dd^cg$ with $g=g_{S,x}\in PSH({\Bbb P}^n,\omega_n)\cap L^\infty_{loc}({\Bbb P}^n\setminus\{x\})$, $g$ has an isotropic pole at $x$ with Lelong number 1 and 
$$(\omega_n+dd^cg)^n=\delta_x.$$
Conversely, any $\omega_n$-psh Green function $g$ with pole at $x$ and maximal Lelong number $\nu(g,x)=1$ is of this form, and in particular it must have an isotropic pole at $x$. Observe that the set of such functions is large. 

\subsubsection{Multipole Green functions}
We push further the result of Proposition \ref{P:mLPn} and study multipole Green functions which arise naturally from rational maps. 

\par Let $f:{\Bbb P}^n\dashrightarrow{\Bbb P}^{n-1}$,
$f=[P_1:\ldots:P_n]$, be a rational map with finite indeterminacy
set $I_f$, where $P_j$ are homogeneous polynomials of degree $d$
on ${\Bbb C}^{n+1}$. Then $f$ determines an $\omega_n$-psh Green
function,
\begin{equation}\label{e:GreenPn}
g_f(\pi_n(z))=d^{-1}\log\|F(z)\|-\log\|z\|,\;z\in{\Bbb
C}^{n+1}\setminus\{0\},\end{equation} where $F:{\Bbb
C}^{n+1}\rightarrow{\Bbb C}^n$, $F(z)=(P_1(z),\dots,P_n(z))$. The
function $g_f$ is continuous, $I_f=\{g_f=-\infty\}$, and $g_f$ has
an isolated pole at each point of $I_f$. Moreover, $g_f$
verifies the Monge-Amp\`ere equation
$$(\omega_n+dd^cg_f)^n=
\sum_{p\in I_f}m_p\delta_p,\;{\rm where}\;m_p>0,\;m_p\in{\Bbb Q},\;\sum_{p\in I_f}m_p=1.$$
Our next result shows that this function has an
extremal property (see \cite{CoMZ} for a similar characterization of classes of pluricomplex Green functions on ${\Bbb C}^n$): 

\begin{Theorem}\label{T:rigid} If $\varphi\in PSH({\Bbb P}^n,\omega_n)$
and $\varphi\leq g_f$, then there exists a unique function $h\in
PSH({\Bbb P}^{n-1},\omega_{n-1})$ such that $\varphi=g_f+d^{-1}h\circ
f$. Conversely, any such function $\varphi$ is $\omega_n$-psh. We have that $\varphi$ is locally bounded on ${\Bbb P}^n\setminus I_f$ if and only if $h$ is bounded. In this case, $\varphi$ satisfies 
$$(\omega_n+dd^c\varphi)^n=\sum_{p\in I_f}m_p\delta_p.$$
\end{Theorem}

\begin{pf} Since the indeterminacy set $I_f$ is finite, we can
find a hyperplane $H$ which does not intersect $I_f$. Let $L$ be a
linear polynomial defining $H$, and let $P_0=L^d$. The map $\hat
f=[P_0:P_1:\ldots:P_n]:{\Bbb P}^n\rightarrow{\Bbb P}^n$ is
holomorphic and $f=\wp\circ\hat f$, where
$$\wp:{\Bbb P}^n\dashrightarrow{\Bbb
P}^{n-1},\;\wp([z_0:z_1:\ldots:z_n])=[z_1:\ldots:z_n],$$ is the
projection with center $[1:0:\ldots:0]$.

\par For every $p\in{\Bbb P}^{n-1}$ the fiber $X_p:=f^{-1}(p)=\hat f^{-1}(\wp^{-1}(p))$ is one-dimensional and is connected by \cite[Proposition 1]{FuH}, since $\wp^{-1}(p)$ is a line in ${\Bbb P}^n$. This implies in particular the uniqueness of $h$.

\par Fix now an arbitrary $p\in{\Bbb P}^{n-1}$, and let us assume
$p=[a_1:\ldots:a_{n-1}:1]$. Then $X_p$ is defined by the equations
$P_j=a_jP_n$. Let $q=[b_0:\ldots:b_n]$ be a point in $X_p\setminus
I_f$. We assume that $b_0=1$. Then $q$ has a neighborhood where
$P_n(1,z_1,\dots,z_n)\neq0$. So, for some constant $c$, we have
$\log\|F\|=\log|P_n|+c$ in this neighborhood. It follows that
$\varphi-g_f$ is psh in some open set which contains $X_p\setminus
I_f$. Since $\varphi-g_f\leq0$ and $I_f$ is a finite set,
$\varphi-g_f$ extends to a subharmonic function on $X_p$. But
$X_p$ is compact and connected, so $\varphi-g_f$ is constant on
$X_p$. We conclude that $\varphi=g_f+(h\circ f)/d$, for some
function $h$ on ${\Bbb P}^{n-1}$. Since $\varphi\leq g_f$ and
$g_f$ is continuous, it follows easily that $h$ is upper
semicontinuous.

\par We now show that $h\in PSH({\Bbb P}^{n-1},\omega_{n-1})$. By using an automorphisms of 
${\Bbb P}^n$ we may assume that the hyperplane $H=\{z_0=0\}$ does not intersect
$I_f$. We claim that the map $F':{\Bbb C}^n\rightarrow{\Bbb C}^n$,
$F'(z')=F(1,z')$, is proper. Indeed, if $P_j^d(z')$ is the
homogeneous part of degree $d$ of $P_j(1,z')$, then $P_j^d(z')$,
$j=1,\dots,n$, have no common zeros except at 0. The homogeneity
of $P_j^d$ yields
$$\sum_{j=1}^n|P_j^d(z')|^2\geq M\|z'\|^{2d},$$ for some constant
$M>0$, which implies that $F'$ is proper. The function
$$u(z')=\varphi([1:z'])+\log\sqrt{1+\|z'\|^2}=
\frac{1}{d}\,\log\|F'(z')\|+\frac{1}{d}\,h\circ\pi_{n-1}(F'(z'))$$
is psh on ${\Bbb C}^n$. Since $F'$ is proper, the function
$$v(w)=d\max\{u(z'):\,F'(z')=w\}=\log\|w\|+h\circ\pi_{n-1}(w)$$
is psh on ${\Bbb C}^n$. This proves that $h\in PSH({\Bbb
P}^{n-1},\omega_{n-1})$.

\par For the converse, note that $$\omega_n+dd^c(g_f+(h\circ
f)/d)=d^{-1}f^\star(\omega_{n-1}+dd^ch)\geq0,$$ so $g_f+(h\circ f)/d$ is
$\omega_n$-psh.

\par Finally, it is clear that $\varphi\in L^\infty_{loc}({\Bbb P}^n\setminus I_f)$ if and only if $h$ is bounded. 
Then we infer by \cite{D93} that
$m_p=(\omega_n+dd^cg_f)^n(\{p\})=(\omega_n+dd^c\varphi)^n(\{p\})$. The
conclusion follows since $\sum_{p\in I_f}m_p=1$.
\end{pf}

\par Note that Proposition \ref{P:mLPn} follows from Theorem \ref{T:rigid} applied to rational maps of degree $d=1$. We will see in section \ref{SS:GP1} that Green functions determined by certain rational maps $f:{\Bbb P}^2\dashrightarrow{\Bbb P}^1$ with three points of indeterminacy provide rich classes of examples of Green functions with one pole on ${\Bbb P}^1\times{\Bbb P}^1$ (see Example \ref{E:GreenP1}). 

\begin{Example}\label{E:homPn} An important particular case of Theorem \ref{T:rigid} is the one of 
rational functions $f:{\Bbb P}^2\dashrightarrow{\Bbb P}^1$, $f=[P_1:P_2]$, where $P_j$ are homogeneous polynomials of degree $d$ whose common zero set $I_f$ consists of $d^2$ distinct points of ${\Bbb P}^2$. Then $g_f$ is a $\omega_2$-psh Green function with $d^2$ isotropic poles and Lelong number $1/d$ 
at each pole. If $d=2$ we observe that any set of four points in general position is the complete intersection of two conics, hence it can be realized as the indeterminacy set $I_f$ for a rational map $f$ of degree $d=2$ as described above. It follows that the $\omega_2$-psh Green functions with four isotropic poles are described by Theorem \ref{T:rigid}. However, if $d\geq3$ a set of $d^2$ points of ${\Bbb P}^2$ in general position is not the complete intersection of two curves of degree $d$ (in fact when $d\geq4$, there is no curve of degree $d$ passing through $d^2$ points in general position). So the Green functions $g_f$ with $d^2$ isotropic poles, $d\geq3$, only exist for very special sets of poles.
\end{Example}

\subsubsection{Partial Green functions}
We compute here in the case of $({\Bbb P}^n,\omega_n)$ the functions $\psi_{\nu,p}$ constructed in Proposition \ref{P:pLelX}. Assume without loss of generality that $p=0\in{\Bbb C}^n$. For $\nu<1$ define $R_\nu,\,C_\nu$ by 
$$R_\nu=[\nu/(1-\nu)]^{1/2} ,\;\nu\log R_\nu+C_\nu=\log\sqrt{1+R_\nu^2}. $$ 
For $z\in{\Bbb C}^n$ let 
$$V(z)=\left\{\begin{array}{ll}
\nu\log\|z\|+C_\nu,\;\|z\|\leq R_\nu,\\
\log\sqrt{1+\|z\|^2}\;,\;\|z\|\geq R_\nu.\end{array}\right.$$

\begin{Proposition}\label{P:compPn} For $\nu<1$ and $z\in{\Bbb C}^n$ we have $\psi_{\nu,p}(z)=V(z)-\log\sqrt{1+\|z\|^2}$.
\end{Proposition}

\begin{pf} Note that $\psi_{\nu,p}(z)=W(z)-\log\sqrt{1+\|z\|^2}$, where 
$$W(z)=\sup\left\{v(z):\,v\in PSH({\Bbb
C}^n),\;v\leq\log\sqrt{1+\|\cdot\|^2},\;\nu(v,0)\geq\nu\right\}.$$
Since $\max_{\|z\|=r}v(z)$ is a convex increasing function of
$\log r$, and since $x=\log R_\nu$ is the solution of the equation
$\frac{d}{dx}\,\log\sqrt{1+e^{2x}}=\nu$, it follows that $W=V$.
\end{pf}

\par Letting $\nu\nearrow1$ it follows that $\psi_{1,p}(z)=\log(\|z\|/\sqrt{1+\|z\|^2})$, $z\in{\Bbb C}^n$, is the Green function constructed in Theorem \ref{T:Kahcon} for $u(z)=\log\|z\|$.

\subsubsection{Dynamical Green functions}\label{SSS:dyn}
We now consider the problem of constructing Green functions on ${\Bbb P}^2$ with one pole at $p$ and Lelong number at $p$ less than 1. Let $\omega=\omega_2$, let $[t:x:y]$ denote the homogeneous coordinates on ${\Bbb P}^2$, and identify $z=(x,y)\in{\Bbb C}^2$ to $[1:x:y]$. Simple examples can be obtained by considering a smooth curve with a flex at $p$, i.e. the tangent line at $p$ does not intersect the curve at any other points. More generally,  for integers $1\leq k<n$, the function 
$$g([t:x:y])=\frac{1}{2n}\log(|y^kt^{n-k}-x^n|^2+|y^n|^2)-\frac{1}{2}\log(|t|^2+|x|^2+|y|^2)$$
is $\omega$-psh and smooth away from $p=0\in{\Bbb C}^2$, $\nu(g,p)=k/n$ and $(\omega+dd^cg)^2=\delta_p$. 

\medskip

\par We describe next more elaborate constructions using complex dynamics. 
Let $h:{\Bbb C}^2\to{\Bbb C}^2$ be a polynomial mapping
of algebraic degree $\lambda>1$. Then $h$ extends to a rational
self-map of ${\Bbb P}^2$, denoted again by $h$, with finite
indeterminacy set $I\subset\{t=0\}$. We call $h$ {\em weakly
regular} if $h$ maps $\{t=0\}\setminus I$ to a point $Z\not\in I$
(see \cite{GSArk}). Such $h$ is algebraically stable ($\deg
h^n=\lambda^n$). It was shown in \cite{S99} that the currents
$\lambda^{-n}(h^n)^\star\omega$ converge weakly to an invariant
positive closed current $T=T_h$ on ${\Bbb P}^2$, $T=\omega+dd^cg$.
We call $T$ the dynamical Green current and $g$ a dynamical Green
function of $h$. By \cite[Theorem 2.2]{GSArk}, $g$ is continuous
on ${\Bbb P}^2\setminus I$, $T\wedge T$ is supported on $I$, so
$g$ is a $\omega$-psh Green function with poles in $I$.

\par If $|I|=1$ then $T\wedge
T=\delta_I$. Our goal is to compute the Lelong number $\nu(T,I)$.

\begin{Proposition}\label{P:dyn} Let $h$ be a weakly regular
polynomial endomorphism of ${\Bbb C}^2$ of degree $\lambda>1$,
with $|I|=1$, and such that
\begin{equation}\label{e:Iest}
dist(h(p),I)\geq C\,dist(p,I)^\delta,\;p\in{\Bbb
P}^2\setminus\{I\},\end{equation} for constants $0<C<1$,
$1<\delta<\lambda$. Then
$\nu(\lambda^{-n}(h^n)^\star\omega,I)\nearrow\nu(T,I)$ as
$n\nearrow\infty$.
\end{Proposition}

\begin{pf} If $\lambda^{-1}h^\star\omega=\omega+dd^c\psi$, where
$\psi\leq0$ is $\omega$-psh, then by \cite[Theorem 2.1]{GAJM}
$$T_n:=\lambda^{-n}(h^n)^\star\omega=\omega+dd^c
g_n\;,\;g_n=\sum_{j=0}^{n-1}\lambda^{-j}\psi\circ h^j\;\searrow\;
g=\sum_{j=0}^\infty\lambda^{-j}\psi\circ h^j,$$ and
$T=\omega+dd^cg$. Hence $\{\nu(T_n,I)\}$ is increasing and
$\nu(T_n,I)\leq\nu(T,I)$.

\par It follows from (\ref{e:Iest}) that there is $C'>0$ so that
for every $n$ and $p\in{\Bbb P}^2\setminus\{I\}$
$$dist(h^n(p),I)\geq(C'\,dist(p,I))^{\delta^n}.$$
Note that the function $\psi$ is smooth except at $I$, and
$\psi\geq\gamma\log\,dist(\cdot,I)-M$ holds on ${\Bbb P}^2$ for
some constants $\gamma,M>0$. Writing $g=g_n+\rho_n$, we deduce that
$$\rho_n(p)\geq\sum_{j=n}^\infty\lambda^{-j}\left(\gamma\log 
dist(h^j(p),I)-M\right)\geq\gamma'(\delta/\lambda)^n\log dist(p,I)
-\epsilon_n,$$ with some $\gamma'>0$ and $\epsilon_n\to0$.
Thus $\nu(T_n,I)\leq\nu(T,I)\leq\nu(T_n,I)+\gamma'(\delta/\lambda)^n$.\end{pf}

\par Note that (\ref{e:Iest}) holds for H\'enon maps $h(x,y)=(P(x)+ay,x)$, $\deg P=\lambda$, 
with $\delta=1$, since $I=[0:0:1]$ is an attracting fixed point for $h^{-1}$. However, the map   
$h(x,y)=(x^\lambda-y^{\lambda-1},y^{\lambda-1})$ shows that (\ref{e:Iest}) does
not hold for $\delta<\lambda$.

\begin{Proposition}\label{P:dyncon} Let $h(x,y)=(x^\lambda+y^\mu,x)$,
where $\lambda>\mu\geq1$ are integers, so $I=[0:0:1]$. The Green
current $T$ of $h$ verifies $T\wedge T=\delta_I$,
$\nu(T,I)=(\lambda-\mu)/\lambda$.\end{Proposition}

\begin{pf} We show first that (\ref{e:Iest}) holds with
$\delta=\lambda-1$. Note that $h$ is weakly regular and in local coordinates $(t,x)$ near $I$ we have
$$h(t,x)=\left(\frac{t}{x}\;,\frac{x^\lambda+
t^{\lambda-\mu}}{xt^{\lambda-1}}\right).$$ It is enough to prove
(\ref{e:Iest}) for $p=(t,x)$ with $0<|x|,|t|<1$. If $|t|\geq|x|$,
or if $|x^\lambda+ t^{\lambda-\mu}|\geq|xt^{\lambda-1}|$, then
$\|h(t,x)\|\geq1$ and the estimate follows. Otherwise, we have
$|t|<|x|<1$ and $|x^\lambda+ t^{\lambda-\mu}|<|xt^{\lambda-1}|$,
so $|x|^\lambda<2|t|^{\lambda-\mu}$. Therefore
$$\|h(t,x)\|\geq\frac{|t|}{|x|}\geq C|x|^{\mu/(\lambda-\mu)}\geq
C|x|^{\lambda-1}\geq C'\,dist(p,I)^{\lambda-1}.$$

\par Next we compute $\nu_n:=\nu(\lambda^{-n}(h^n)^\star\omega,I)$. Let  $h^n([t:x:y]=[t^{\lambda^n}:p_n(t,x,y):q_n(t,x,y)]$,
where $p_n,q_n$ are homogeneous polynomials of degree $\lambda^n$, and
$$v_n(t,x)=\log(
|t|^{2\lambda^n}+|p_n(t,x,1)|^2+|q_n(t,x,1)|^2)^{1/2}\,.$$ 
It follows by induction that
$\nu(v_n,0)=\lambda^n-\max\{\deg_yp_n,\deg_yq_n\}=\lambda^n-
\mu\lambda^{n-1}$, where $\deg_yp_n$ denotes the degree in $y$ of
$p_n$. Hence $\nu_n=(\lambda-\mu)/\lambda=\nu(T,I)$.\end{pf}

\par If $h$ is H\'enon map of degree $\lambda$ a similar argument shows 
$\nu(T_h,I)=1-\lambda^{-1}$.

\section{Green functions on ${\Bbb P}^1\times{\Bbb P}^1$}
\label{S:P1P1}
\par It is possible to describe the functions $\nu,\,\varepsilon,\,g,\,\psi$ on a multiprojective 
space ${\Bbb P}^{n_1}\times\dots\times{\Bbb P}^{n_k}$. For simplicity, 
we only consider the case $X={\Bbb P}^1\times{\Bbb P}^1={\Bbb P}^1_z\times{\Bbb P}^1_w$. Let $\pi_z:X\to{\Bbb P}^1_z$, 
$\pi_w:X\to{\Bbb P}^1_w$, denote the canonical projections and set 
$$\alpha_{a,b}:=a\alpha_z+b\alpha_w,\;\omega_{a,b}:=a\omega_z+b\omega_w,\;a,b\geq0,$$
where $\alpha_z=\pi_z^\star\alpha_1$, $\alpha_w=\pi_w^\star\alpha_1$, $\omega_z=\pi_z^\star\omega_1$, 
$\omega_w=\pi_w^\star\omega_1$, and $\omega_1\in\alpha_1$ is the Fubini-Study form on ${\Bbb P}^1$. 
Note that $\alpha_{a,b}$ is a K\"ahler class if and only if $a,b>0$. 

\par For concrete computations, it will be convenient to use coordinates on $X$. Let
$$\pi:({\Bbb C}^2\setminus\{0\})\times({\Bbb C}^2\setminus\{0\})\rightarrow
X\;,\;\pi(z_0,z_1,w_0,w_1)=([z_0:z_1],[w_0:w_1]),$$ 
and identify $(z_1,w_1)\in{\Bbb C}^2$ to $\pi(1,z_1,1,w_1)\in X$. The currents $T\in{\mathcal P}(\alpha_{a,b})$ can be described using the class $P_{a,b}$ of bihomogeneous psh
functions $\widetilde u$ on ${\Bbb C}^4$ (see \cite{GAJM}): 
$$\widetilde u(\lambda z_0,\lambda z_1,\mu w_0,\mu w_1)=a\log|\lambda|+b\log|\mu|+\widetilde
u(z_0,z_1,w_0,w_1),\;\lambda,\mu\in{\Bbb C}.$$ 
Then $\pi^\star T=dd^c\widetilde u$, for some $\widetilde u\in P_{a,b}$ which is unique up to additive constants.  

\par For a point $p=(x,y)\in X$ we denote by
$$V_x=\pi_z^{-1}(x)=\{z=x\}\;,\;H_y=\pi_w^{-1}(y)=\{w=y\},$$ 
the vertical, and respectively horizontal, line through $p$.

\subsection{Maximal Lelong numbers}\label{SS:P1Max}

\begin{Proposition}\label{P:LelP1} For all $p=(x,y)\in X$, we have 
$$\nu(\alpha_{a,b},p)=a+b,\;\varepsilon(\alpha_{a,b},p)=\min\{a,b\}.$$
If $T\in{\mathcal P}(\alpha_{a,b})$ and $\nu(T,p)=a+b$ then $T=a[V_x]+b[H_y]$. Moreover, if $T$ does not charge $V_x$ and $H_y$ then $\nu(T,p)\leq\min\{a,b\}$.
\end{Proposition}

\begin{pf} Let $T\in{\mathcal P}(\alpha_{a,b})$. We can assume that $p=(0,0)$ and let $m=\min\{a,b\}$. The current $R_{a,b}\in{\mathcal P}(\alpha_{a,b})$ 
defined by $\pi^\star R_{a,b}=dd^c\widetilde u_{a,b}$, where $\widetilde u_{a,b}\in P_{a,b}$,
$$\widetilde u_{a,b}(z_0,z_1,w_0,w_1):=m\log\sqrt{|z_1w_0|^2+|w_1z_0|^2}+(a-m)\log|z_0|+
(b-m)\log|w_0|,$$
shows that $\varepsilon(\alpha_{a,b},p)\geq m$. Moreover, the measure $T\wedge R_{1,1}$ is well
defined and
$$\nu(T,p)=T\wedge R_{1,1}(\{p\})\leq\int_XT\wedge R_{1,1}=\int_X\omega_{a,b}\wedge\omega_{1,1}=a+b.$$

\par Assume now that $T$ does not charge the subvarieties $V_x$
and $H_y$. By \cite{D92}, there exist $\epsilon_j\searrow0$ and currents 
$T_j\in{\mathcal P}(\alpha_{a,b}+\epsilon_j\alpha_{1,1})$ with analytic singularities, so that $0\leq\nu(T,q)-\nu(T_j,q)\le\epsilon_j$ for every $q\in X$. Since $T$ does not charge $V_x$, the measure $T_j\wedge[V_x]$ is 
well defined. If $v_j$ is a psh potential of $T_j$ near $p$ then
$$\nu(T_j,p)\leq\nu(v_j\vert_{_{V_x}},p)=
T_j\wedge[V_x](\{p\})\leq\int_XT_j\wedge
[V_x]=b+\epsilon_j.$$ We replace $V_x$ by $H_y$ in this argument and let $j\to+\infty$ to get  
$\nu(T,p)\leq m$. By (\ref{e:Seshadri}) it follows that 
$\varepsilon(\alpha_{a,b},p)\leq m$. 

\par Assume finally that $\nu(T,p)=a+b$. By \cite{Siu74}, we can write
$$T=a'[V_x]+b'[H_y]+T',\;T'\in{\mathcal
P}(\alpha_{a-a',b-b'}),$$ where $T'$ does not charge $V_x$ and $H_y$. By what we have already shown, 
$$a+b=\nu(T,p)\leq a'+b'+\min\{a-a',b-b'\}.$$ This implies that
$a'=a$, $b'=b$, and $T'=0$. 
\end{pf}

\par Observe that the functions $\nu,\,\varepsilon$ are constant here, as well as in the case of ${\Bbb P}^n$, because $\alpha$ is invariant under a compact group of automorphisms that acts transitively on $X$. 

\par Note that ${\rm Vol}\,(\alpha_{a,b})^{1/2}=\sqrt{2ab}>\min\{a,b\}$, hence the upper bound given in (\ref{e:Seshvar}) is not sharp in this case. Another obvious consequence of the previous proposition is the following: 

\begin{Corollary}\label{C:LelP1} There is no Green function with one isotropic pole on ${\Bbb P}^1\times{\Bbb P}^1$.
\end{Corollary}

\par We can however compute the partial Green functions with isotropic singularity $\psi_{\nu,p}$ constructed in Proposition \ref{P:pLelX}. Assume that $p=(0,0)\in{\Bbb C}^2\subset X$, and let 
$a=b=1$, $\nu=\varepsilon(\alpha_{1,1},p)=1$. A psh potential of $\omega_{1,1}$ on 
${\Bbb C}^2$ is given by  
$$\rho(z_1,w_1)=\log\sqrt{1+|z_1|^2}+\log\sqrt{1+|w_1|^2}.$$

\begin{Proposition}\label{P:compP1} We have $\psi_{1,p}(z_1,w_1)=
\log(|z_1|+|w_1|)-\rho(z_1,w_1)$ if $|z_1w_1|\leq1$, and $\psi_{1,p}(z_1,w_1)=0$ if $|z_1w_1|\geq1$.
\end{Proposition}

\begin{pf} We have to obtain upper estimates for psh functions $v$
on ${\Bbb C}^2$ which verify $v\leq\rho$ and $\nu(v,0)\geq1$. We
do this first along a complex line $z_1=s\zeta$, $w_1=t\zeta$.
Using the same convexity argument as in the proof of Proposition \ref{P:compPn},
we obtain
$$v(s\zeta,t\zeta)\leq\left\{\begin{array}{ll}
\log|\zeta|+C,\;|\zeta|\leq R,\\
\rho(s\zeta,t\zeta),\;|\zeta|\geq R.\end{array}\right.$$ Here
$R=|st|^{-1/2}$, $x=\log R$ is the solution of the equation
$$\frac{d}{dx}\left(\log\sqrt{1+|s|^2e^{2x}}+
\log\sqrt{1+|t|^2e^{2x}}\right)=1,$$ and $C=\log(|s|+|t|)$
verifies $\log R+C=\rho(sR,tR)$. If $s=1,\,t=w_1/z_1$, we get   
$$v(z_1,w_1)\leq V(z_1,w_1)=\left\{\begin{array}{ll}
\log(|z_1|+|w_1|),\;|z_1w_1|\leq1,\\
\rho(z_1,w_1),\;|z_1w_1|\geq1.\end{array}\right.$$ Since
$\log(|z_1|+|w_1|)\leq\rho(z_1,w_1)$ on ${\Bbb C}^2$, with
equality when $|z_1w_1|=1$, the function $V$ is psh. It follows
that $\psi_{1,p}=V-\rho$.  
\end{pf}

\par Note that the (unbounded) hyperconvex domain 
$$D_{1,p}=\{\psi_{1,p}<0\}=\{(z_1,w_1)\in{\Bbb C}^2:\,|z_1w_1|<1\}$$ 
does not have a pluricomplex Green function: if $v<0$ is psh on $D_{1,p}$ and $v(0,0)=-\infty$ 
then $v=-\infty$ along the lines $\{z_1=0\},\{w_1=0\}$.

\subsection{Green functions with one pole}\label{SS:GP1}
It is clear from Proposition \ref{P:LelP1} and Corollary
\ref{C:LelP1} that the characterization of Green functions
in $PSH(X,\omega_{a,b})$ with one pole at $p\in X$ is more
involved. Using a birational map, we will show that they correspond to a certain class of Green functions with three poles on ${\Bbb P}^2$.  A rich class of examples of the latter can be constructed using (\ref{e:GreenPn}) (see also Theorem \ref{T:rigid}). This will show that the Green functions of $X$ with pole 
at $p$ have many different types of singularities, even if one asks that the Lelong number at $p$ is maximal.

\smallskip

We may assume that $p=(0,0)\in{\Bbb C}^2\subset 
X$ and $a=1\leq b$. Let $\omega=\omega_{FS}$ on ${\Bbb P}^2$ and consider the rational map
$\Phi:{\Bbb P}^2\dashrightarrow X$ defined by
$$
\Phi([t_0:t_1:t_2])=([t_0:t_1],[t_0:t_2]).
$$
It is a birational map, with rational inverse
$$
\Phi^{-1}([z_0:z_1],[w_0:w_1])=[z_0w_0:z_1w_0:w_1z_0].
$$ 
Note that
$\Phi$ is the identity on ${\Bbb C}^2\equiv\{[1:t_1:t_2]\in{\Bbb
P}^2\}\equiv\{([1:z_1],[1:w_1])\in X\}$, $\Phi$ blows up the
points $A=[0:1:0]$, $B=[0:0:1]$, to the lines $\{z_0=0\}$,
respectively $\{w_0=0\}$, and $\Phi$ contracts the line $\{t_0=0\}$ to the point $q=(\infty,\infty)$.

\smallskip

\par We denote by ${\mathcal S}_b$ the set of the currents $S\in{\mathcal
P}(\alpha_{1,b})$ with locally bounded potentials on $X\setminus\{p\}$ and such that $S\wedge S=2b\delta_p$.  A potential of $S$ is then a $\omega_{1,b}$-psh Green function on $X$ with pole at $p$.

\par Let ${\mathcal R}_b$ be the set of currents $R\in{\mathcal P}((1+b)\omega)$ on ${\Bbb P}^2$ whose potentials are locally bounded on ${\Bbb P}^2\setminus\{p,A,B\}$, have isotropic poles at $A,B$ with Lelong numbers $\nu(R,A)=b$, $\nu(R,B)=1$, and such that $R\wedge R=0$ on ${\Bbb P}^2\setminus\{p,A,B\}$. It follows that a potential $v$ of $R$ is a $(1+b)\omega$-psh Green function on ${\Bbb P}^2$ with poles at $p,A,B$: 
$$R\wedge R=((1+b)\omega+dd^cv)^2=b^2\delta_A+\delta_B+2b\,\delta_p.$$

\begin{Proposition}\label{P:GreenP1} 
The mapping $\Phi^\star:{\mathcal
S}_b\rightarrow{\mathcal R}_b$ is well defined and bijective. 
Its inverse is the mapping 
$$
G:R \in {\mathcal R}_b \mapsto (\Phi^{-1})^\star R-b[z_0=0]-[w_0=0] \in {\mathcal S}_b.
$$ 
\end{Proposition}

\begin{pf} Let $S\in{\mathcal S}_b$ and $\widetilde u\in P_{1,b}$ be a potential of $S$. Then 
$$\widetilde v(t_0,t_1,t_2):=\widetilde u(t_0,t_1,t_0,t_2),\;\widetilde v(\lambda
t_0,\lambda t_1,\lambda t_2)=\widetilde v(t_0,t_1,t_2)+(1+b)\log|\lambda|,$$
is a logarithmically homogeneous potential for $R=\Phi^\star S$, so $R\in{\mathcal P}((1+b)\omega)$. In particular, it follows that $R$ has locally bounded potentials on ${\Bbb P}^2\setminus\{p,A,B\}$. Near the point $A$, assuming wlog that $|t_0|\leq|t_2|$ we have 
$$\widetilde v(t_0,1,t_2)=\widetilde u(t_0,1,t_0/t_2,1)+b\log|t_2|=b\log\sqrt{|t_0|^2+|t_2|^2}+O(1).$$
So $R$ has potentials with an isotropic pole at $A$ and $\nu(R,A)=b$. One proves in the same way that $R$ has potentials with an isotropic pole at $B$ and $\nu(R,B)=1$. We have $R\wedge R=S\wedge S=0$ on ${\Bbb C}^2\setminus\{0\}$. Since $R$ has locally bounded potentials near each point of $\{t=0\}\setminus\{A,B\}$ we have $R\wedge R(\{t=0\}\setminus\{A,B\})=0$, so $R\in{\mathcal R}_b$.

\par Conversely, let $R\in{\mathcal R}_b$ with logarithmically homogeneous potential $\widetilde v$. Then 
$$\widetilde u(z_0,z_1,w_0,w_1):=\widetilde v(z_0w_0,z_1w_0,w_1z_0)-b\log|z_0|-\log|w_0|\in P_{1,b}$$ 
is a bihomogeneous potential of $G(R)$. We show that $G(R)$ has locally bounded potentials in a neighborhood of any point at infinity $\zeta\neq q$. Suppose wlog $\zeta\in\{z_0=0\}$. Then for $|z_0|$ small enough we have that $[z_0:1:z_0w_1]$ is near $A$, so  
$$\widetilde u(z_0,1,1,w_1)=\widetilde v(z_0,1,w_1z_0)-b\log|z_0|=b\log\sqrt{1+|w_1|^2}+O(1)=O(1).$$
\par Next we study the potentials of $G(R)$ in a neighborhood of $q$. We have 
$$\widetilde u(z_0,1,w_0,1)=\widetilde v(z_0w_0,w_0,z_0)-b\log|z_0|-\log|w_0|,$$
where $|z_0|,|w_0|$ are small. If $|w_0/z_0|$ is small, then $[w_0:w_0/z_0:1]$ is near $B$ so 
$$\widetilde u(z_0,1,w_0,1)=\widetilde v(w_0,w_0/z_0,1)+\log|z_0|-\log|w_0|=\log\sqrt{|z_0|^2+1}+O(1).$$
Similarly, $\widetilde u(z_0,1,w_0,1)=O(1)$ if $|z_0/w_0|$ is small. If $\epsilon\leq|w_0/z_0|\leq M$ then 
$$\widetilde u(z_0,1,w_0,1)=\widetilde v(w_0,w_0/z_0,1)+\log(|z_0|/|w_0|)=O(1).$$ 
It follows that $G(R)$ has locally bounded potentials in $X\setminus\{p\}$, hence $G(R)\in{\mathcal S}_b$.

\par Since $\Phi$ is the identity on ${\Bbb C}^2$ and the currents in ${\mathcal R}_b$, resp. ${\mathcal S}_b$, do not charge the line(s) at infinity, we conclude by the support theorem that $\Phi^\star$ is bijective and $G$ is its inverse. \end{pf}

\begin{Example}\label{E:GreenP1} Let $1\leq b=m/n\in\Bbb Q$ and $f=[P_1:P_2]:{\Bbb P}^2\dashrightarrow{\Bbb P}^1$, where $$P_1(t_0,t_1,t_2)=t_1^{nk}t_2^{mk},\;P_2(t_0,t_1,t_2)=t_1^{nk}t_0^{mk}+t_2^{mk}t_0^{nk}+t_1t_2Q(t_0,t_1,t_2),$$ $k\geq1$ is an integer, and $Q$ is a homogenous polynomial of degree $(m+n)k-2$ with $\deg_{t_1}Q\leq nk-1$ and $\deg_{t_2}Q\leq mk-1$. Note that the indeterminacy set $I_f=\{p,A,B\}$ and the current 
$$R_f:=(1+b)(\omega+dd^cg_f)\in{\mathcal R}_b,$$ where $g_f$ is the Green function associated to $f$ defined in (\ref{e:GreenPn}). Then $S_f=G(R_f)$ has bihomogeneous potential $\widetilde u_f\in P_{1,b}$ given by 
$$\widetilde u_f(1,z_1,1,w_1)=
\frac{1}{2nk}\log\left(|z_1^{nk}w_1^{mk}|^2+|z_1^{nk}+w_1^{mk}+z_1w_1Q(1,z_1,w_1)|^2\right),$$ 
where 
$Q(1,z_1,w_1)=\sum_{i_1=0}^{nk-1}\sum_{i_2=0}^{mk-1}c_{i_1i_2}z_1^{i_1}w_1^{i_2}$.
Depending on the vanishing order of $Q(1,\cdot)$ at the origin, one sees that the Lelong number $\nu(S_f,p)$ can take any value of the form $\frac{j}{nk}$, $2\leq j\leq nk$. It follows that for any rational number $r\in(0,1]$ there exist  $\omega_{1,b}$-psh Green functions on $X$ with one pole at $p$ and Lelong number equal to $r$ there, but with different types of singularities at $p$.
\end{Example}

\par We finally give an alternate way to construct $\omega_{1,1}$-psh Green functions on $X$ with pole at $q=(\infty,\infty)$, using currents on ${\Bbb P}^2$ arising from psh functions in the Lelong class ${\mathcal L}^\star({\Bbb C}^2)$. This is the class of psh functions $v$ on
${\Bbb C}^2$ so that
$$\limsup_{\|s\|\rightarrow\infty}v(s)/\log\|s\|=1.$$ If $R$ is the trivial extension of $dd^cv$ to ${\Bbb P}^2$ then  $R\in{\mathcal P}({\omega})$. 

\begin{Proposition}\label{P:GreenP1a} Let $R\in{\mathcal P}(\omega)$ be a current with locally bounded potentials in ${\Bbb P}^2\setminus\{t_0=0\}$ and near the points $A,B$. Then the current $S=(\Phi^{-1})^\star R\in{\mathcal P}(\alpha_{1,1})$, $\nu(S,q)=1$, and $S$ has  locally bounded potentials on $X\setminus\{q\}$. Moreover, we have  
$$S\wedge S=2\delta_q \Longleftrightarrow R\wedge R=0\;{\rm on}\;{\Bbb
P}^2\setminus\{t_0=0\}.$$
\end{Proposition}

\begin{pf} By considering (bi)homogeneous potentials as in the proof of Proposition \ref{P:GreenP1}, it follows that $S\in{\mathcal P}(\alpha_{1,1})$ and $S$ has  locally bounded potentials on $X\setminus\{q\}$. So  
$S\wedge S(\{z_0=0\}\cup\{w_0=0\}\setminus\{q\})=0$, and $S\wedge S=0$ on ${\Bbb C}^2$ implies $S\wedge S=2\delta_q$.

\par Let $\nu:=\nu(S,q)$. Since $\Phi$ contracts the line $\{t_0=0\}$ to $q$, we have that $\Phi^\star
S=\nu[t_0=0]+T$, where $T\in{\mathcal P}((2-\nu)\omega)$ does not charge the line
$\{t_0=0\}$. Note that $R=T$ on ${\Bbb C}^2$. By the support
theorem we conclude that $R=T$, so $\nu=1$. \end{pf}

\par Proposition \ref{P:GreenP1a} shows how Green
functions can be constructed on $X$ by using 
currents $R$ on ${\Bbb P}^2$ possessing the right properties at
any two points $A,\,B$ and outside the line joining them. Indeed, we pull back $R$ by an automorphism of ${\Bbb P}^2$ which maps the points $[0:1:0],\,[0:0:1]$ to $A,B$, and then apply Proposition \ref{P:GreenP1a}.

\begin{Example}\label{E:GreenP1a} The Green currents
$T^+,T^-$ of a H\'enon map $h$ on ${\Bbb C}^2$ yield by the
preceding considerations Green functions on $X$ with pole at $q$.
More generally, let $h$ be a weakly regular polynomial
endomorphism of ${\Bbb C}^2$ with indeterminacy set $I$ (see
section \ref{SSS:dyn}). Then its Green current $T$ has continuous
local potentials on ${\Bbb P}^2\setminus I$ and $T\wedge
T=\sum_{s\in I}m_s\delta_s$. So $T$ yields a Green function on $X$
with pole at $q$. 
\end{Example}

\section{Del Pezzo Surfaces}\label{S:DelPezzo}
We evaluate here the functions $\nu,\,\varepsilon,\,g$ when $X$ is a (smooth) Del Pezzo surface, 
i.e. $\dim_{\Bbb C}X=2$ and $c_1(X)>0$. It is well known (see e.g. \cite{Dem}) that such $X$ is biholomorphic to either ${\Bbb P}^1\times{\Bbb P}^1$, ${\Bbb P}^2$, or ${\Bbb P}^2$ blown up at $r$ points in general position, $1\leq r\leq8$. Here general position means the following:
\par $-$ no three points are collinear;
\par $-$ no six points lie on a conic;
\par $-$ when $r=8$, the points do not lie on a cubic that is singular at one of them.

\smallskip

\par The cases $X={\Bbb P}^2$, $X={\Bbb P}^1\times{\Bbb P}^1$, have already been considered in Sections \ref{S:Pn} and \ref{S:P1P1}. We focus here on the case when $X$ is the blow up of ${\Bbb P}^2$ at 8 points in general position, which we consider to be the most interesting one. The other cases could be handled similarly. Note that the Seshadri constants $\varepsilon$ are computed in \cite{Br}.

\subsection{Maximal Lelong numbers}\label{SS:DPMax}
Let $\pi:X\to{\Bbb P}^2$ be the blow up of ${\Bbb P}^2$ at 8 points $p_1,\dots,p_8$ in general position, and let $E_j=\pi^{-1}(p_j)$ denote the exceptional divisors. We let 
$$\alpha:=c_1(X)=K_X^{-1}=\pi^\star{\mathcal O}(3)-\sum_{j=1}^8E_j$$
denote the (ample) anticanonical class of $X$. It is well known \cite{Dem} that $2\alpha$ is very ample. It follows from Remark \ref{R:vamp} that 
\begin{equation}\label{e:lowbd}
\nu(\alpha,x)\geq1,\;\varepsilon(\alpha,x)\geq1/2,\;\forall\,x\in X.
\end{equation}

\par We can actually be much more precise. Let ${\mathcal V}$ be the pencil of cubics in ${\Bbb P}^2$ passing through $p_1,\dots,p_8$. It contains at most 12 singular cubics \cite{Dem}. We let $S\subset X$ denote the set of the corresponding singular points, $|S|\leq12$. These points do not belong to the exceptional divisors, by the general position assumption. 

\begin{Proposition}\label{P:LX8} We have 
$$\nu(\alpha,x)=\left\{\begin{array}{ll} 
1,\;\mbox{if $x\in X\setminus S$},\\ 2,\;\mbox{if $x\in S$}. 
\end{array}\right.$$
Moreover, if $x\in S$ and $T\in{\mathcal P}(\alpha)$ does not charge the strict transform of the singular cubic in $\mathcal V$ passing through $x$ then $\nu(T,x)\leq1/2$.
\end{Proposition}

\begin{pf} For $x\in X$ there exists a unique cubic ${\mathcal C}_x\in{\mathcal V}$ whose strict transform ${\mathcal C}'_x$ contains $x$. (If $x\in E_j$ this is the cubic whose strict transform intersects $E_j$ at $x$.) Note that ${\mathcal C}'_x$ is irreducible. 

\par Let $T\in{\mathcal P}(\alpha)$. We assume at first that $T$ does not charge ${\mathcal C}'_x$ and let 
$\omega$ be a fixed K\"ahler form on $X$. By \cite{D92} there exist $\epsilon_j\searrow0$ and currents 
$T_j\in{\mathcal P}(\alpha+\epsilon_j\omega)$ 
with analytic singularities, such that $T_j\to T$ and $0\leq\nu(T,z)-\nu(T_j,z)\leq\epsilon_j$ for all $z\in X$. 
Since $T$ does not charge ${\mathcal C}'_x$, the measure $T_j\wedge[{\mathcal C}'_x]$ is 
well defined. As ${\rm Vol}\,(\alpha)=1$ it follows that 
$$1+O(\epsilon_j)=\int_X T_j\wedge[{\mathcal C}'_x]
\geq T_j\wedge[{\mathcal C}'_x](\{x\})\geq\nu(T_j,x)m({\mathcal C}'_x,x),$$
where $m({\mathcal C}'_x,x)$ denotes the multiplicity of ${\mathcal C}'_x$ at $x$. 
The last inequality can be seen by using a local normalization at $x$ for each irreducible component of  
${\mathcal C}'_x$ and since local psh potentials of $T_j$ are subharmonic along ${\mathcal C}'_x$. 

\par Letting $j\to+\infty$, we have shown that $\nu(T,x)\leq1/m({\mathcal C}'_x,x)\leq1$, if $T\in{\mathcal P}(\alpha)$ does not charge ${\mathcal C}'_x$. In particular, if $x\in S$ then $\nu(T,x)\leq1/2$ since $m({\mathcal C}'_x,x)=2$.

\par In the general case, we can write by \cite{Siu74}
$$T=a[{\mathcal C}'_x]+(1-a)R,\;0\leq a\leq1,$$ where $R\in{\mathcal P}(\alpha)$ does not charge 
${\mathcal C}'_x$. Then 
$$\nu(T,x)=am({\mathcal C}'_x,x)+(1-a)\nu(R,x)\leq a(m({\mathcal C}'_x,x)-1)+1\leq 
m({\mathcal C}'_x,x),$$ 
which concludes the proof. 
\end{pf}

\subsection{Uniform integrability exponent}\label{SS:DPKE}
We fix $\omega\in\alpha=c_1(X)$ a K\"ahler form and we denote by $PSH_0(X,\omega)$ the set of $\omega$-psh functions $\varphi$ normalized by $\max_X\varphi=0$. This is a compact subset of $L^1(X)$. Set 
$$\sigma(X)=\sup\{c\geq0:\,e^{-2c\varphi}\in L^1(X),\;\forall\,\varphi\in PSH_0(X,\omega)\}.$$
This number clearly depends only on $\alpha=c_1(X)$, rather than on the particular choice of $\omega$. By the compactness of $PSH_0(X,\omega)$ and the semicontinuity of the ``complex singularity exponent" \cite{DK01}, $\sigma(X)$ coincides with the exponent introduced by Tian in \cite{Ti87} (the so-called ``$\alpha$-invariant of Tian"). 

\smallskip

We assume here again that $X$ is the blow up of ${\Bbb P}^2$ at $8$ points in general position.
Since $\nu(\alpha,x)\leq2$ for all $x\in X$, it follows from Skoda's integrability theorem \cite{Sk72} that 
$\sigma(X)\geq1/2$. One can however obtain sharp estimates, thanks to the full characterization given in Proposition \ref{P:LX8}:

\begin{Proposition}\label{P:TX8} If there is a singular cubic in ${\mathcal V}$ with a cusp then $\sigma(X)=5/6$. Otherwise, $\sigma(X)=1$. 
\end{Proposition}

\par Recall that there is no cuspidal cubic in ${\mathcal V}$ when the points $p_1,\dots,p_8$ are in very general position \cite{Dem}.

\medskip

\par\noindent{\em Proof of Proposition \ref{P:TX8}.} Let $s=|S|\leq12$ and ${\mathcal C}'_j$, $1\leq j\leq s$, denote the strict transforms of the singular cubics in ${\mathcal V}$. We write $[{\mathcal C}'_j]=\omega+dd^c\varphi_j$, where $\varphi_j\in PSH_0(X,\omega)$. 

\par Fix now $\varphi\in PSH_0(X,\omega)$ and let $T=\omega+dd^c\varphi\in{\mathcal P}(\alpha)$. By \cite{Siu74}, 
$$T=a_0T_0+\sum_{j=1}^sa_j[{\mathcal C}'_j],\;{\rm where}\;a_j\geq0,\;\sum_{j=0}^sa_j=1,$$
and $T_0=\omega+dd^c\varphi_0\in{\mathcal P}(\alpha)$ does not charge any curve ${\mathcal C}'_j$. 
H\"older's inequality shows that $e^{-2c\varphi}\in L^1(X)$ if $e^{-2c\varphi_j}\in L^1(X)$ for all $j=0,\dots,s$. 

\par For $j\geq1$, a direct computation in local coordinates shows that $e^{-2c\varphi_j}\in L^1(X)$ for every $c<1$ if ${\mathcal C}'_j$ is non-singular or has a simple node, while $e^{-2c\varphi_j}\in L^1(X)$ for every $c<5/6$ if ${\mathcal C}'_j$ has a cusp. In the latter case, $e^{-2c\varphi_j}\not\in L^1(X)$ if $c=5/6$.

\par Since $T_0$ does not charge any curve ${\mathcal C}'_j$, it follows from Proposition \ref{P:LX8} that 
$\nu(T_0,x)\leq1$ for all $x\in X$. By \cite{Sk72} we see that  $e^{-2c\varphi_0}\in L^1(X)$ for every $c<1$. This completes the proof of the proposition. $\Box$

\medskip

\par Note that $\sigma(X)$ is also called the (global) ``log-canonical threshold" of $X$. 
It has been the subject of intensive studies in the last decade. 
The above result has been recently obtained by Cheltsov \cite{Che} by more algebraic methods. 

\smallskip

The importance of this notion is seen in its connection with the existence
of K\"ahler-Einstein metrics: it was shown by Tian \cite{Ti87}
that a Fano surface admits a K\"ahler-Einstein metric if $\sigma(X)>2/3$.
The exponent $\sigma(X)$ was previously estimated by Tian and Yau
in \cite{TY}.

\subsection{Green functions}\label{SS:DPGreen}

In this section $X$ denotes again the blow up of ${\Bbb P}^2$ at $8$ points in general position.

\subsubsection{Special points} 
For $x\in S$, let ${\mathcal C}_x$ be the cubic in ${\mathcal V}$ which is singular at $x$, and 
let ${\mathcal C}'_x$ be its strict transform. 

\par Counting dimension we see that there exists an irreducible sextic $Z\subset{\Bbb P}^2$ passing through 
$x$ and with multiplicity 2 at each point $p_j$. By Bezout we see that $Z$ and ${\mathcal C}_x$ intersect only at $x$ and at the points $p_j$ and the intersection numbers $(Z\cdot{\mathcal C}_x)_{p_j}=(Z\cdot{\mathcal C}_x)_x=2$. This implies that the strict transform $Z'\subset X$ of $Z$ intersects ${\mathcal C}'_x$ only at $x$ with $(Z'\cdot{\mathcal C}'_x)_x=2$. 

\par We write $(1/2)[Z']=\omega+dd^cu,\;[{\mathcal C}'_x]=\omega+dd^cv$, and set 
$$g_x:=(1/2)\log(e^{2u}+e^{2v})\in PSH(X,\omega)\cap C^\infty(X\setminus\{x\}).$$

\begin{Proposition}\label{P:GX8} If $x\in S$ we have $(\omega+dd^cg_x)^2=\delta_x$, and the function $g_x$ is a $\omega$-psh Green function with Lelong number $\nu(g_x,x)=1/2$.
\end{Proposition}

\begin{pf} Since $Z'$ is smooth at $x$ we have $\nu(g_x,x)=1/2$. Moreover, $(Z'\cdot{\mathcal C}'_x)_x=2$ implies that $(\omega+dd^cg_x)^2(\{x\})=1$. We conclude by mass considerations.
\end{pf}

\par Observe that the singularity of $g_x$ at $x$ is not isotropic, since an isotropic pole with Lelong number $1/2$ would produce a Dirac mass at $x$ with coefficient $1/4$. However, the existence of a Green function which is locally bounded away from $x$ has interesting consequences:

\begin{Corollary}\label{C:GX8} If $x\in S$ then $\varepsilon(\alpha,x)=1/2$. Moreover, the 
supremum is attained in the formula (\ref{e:Seshadri}) of $\varepsilon(\alpha,x)$, i.e. 
$$\exists\,\varphi\in PSH(X,\omega)\cap L^\infty_{loc}(X\setminus\{x\}),\;\|\varphi-(1/2)\log dist(\cdot,x)\|_{L^\infty(X)}<+\infty.$$
\end{Corollary}

\begin{pf} It follows from (\ref{e:lowbd}) and Proposition \ref{P:LX8} that 
$\varepsilon(\alpha,x)=1/2$. Let $g_x$ be the function constructed in Proposition \ref{P:GX8}. Fix $\chi\in C^\infty(X)$ a test function with $\chi\equiv 1$ on $\overline U$, where $U$ is a small open neighborhood of $x$. We define 
$$\varphi:=\max\{g_x,(1/2)\chi\log dist(\cdot,x)-C\},$$ 
where $C$ is large so that $\varphi=g_x$ on $X\setminus U$. Since $\chi\log dist(\cdot,x)$ is psh on $U$ we see that $\varphi\in PSH(X,\omega)$. Now $\nu(g_x,x)=1/2$, therefore $\varphi-(1/2)\log dist(\cdot,x)$ is bounded on $X$.
\end{pf}

\subsubsection{Generic points}\label{SSS:DPgen}
Assume now that $x\in X\setminus S$. The bound (\ref{e:lowbd}) is not sharp: by \cite{Br} we have 
$\varepsilon(\alpha,x)=1$. 

\par It is easy to see that the supremum in formula (\ref{e:Seshadri}) is attained if $x$ is the ninth base point of the pencil of cubics ${\mathcal V}$. In this case we write $[{\mathcal C}'_1]=\omega+dd^cu,\;[{\mathcal C}'_2]=\omega+dd^cv$, where ${\mathcal C}'_j$ are the strict transforms of two cubics generating ${\mathcal V}$, and we set 
$$g_x:=(1/2)\log(e^{2u}+e^{2v})\in PSH(X,\omega)\cap C^\infty(X\setminus\{x\}).$$
We have that $(\omega+dd^cg_x)^2=\delta_x$ and $g_x$ is a $\omega$-psh Green function with an isotropic pole at $x$ with $\nu(g_x,x)=1$.

\par However, it is unclear whether this holds at arbitrary points $x\in X\setminus S$. If this was the case, it would imply that $K_Y^{-1}$ admits a positive metric with bounded potentials, where $Y\to{\Bbb P}^2$ is the blow up of ${\Bbb P}^2$ at 9 points in general position, which is a famous open problem (see \cite{DPS}). 
Observe that the existence of such a metric is equivalent to constructing a 
$\omega_{FS}$-psh Green function with isotropic poles of Lelong number $1/3$ at 9 points in general position 
in ${\Bbb P}^2$.

\par More generally, finding a $\omega_{FS}$-psh Green function with isotropic poles of Lelong number 
$1/\sqrt{s}$ at $s$ points in general position in ${\Bbb P}^2$ is equivalent to the celebrated (strong version of) Nagata's conjecture (see \cite[Remark 5.1.14]{La}).

\end{document}